\documentclass[11pt]{amsart}
\usepackage{amsfonts}
\usepackage{amsmath}
\usepackage{amssymb}
\usepackage{amsthm}
\usepackage{ color, ulem}

\input xy
\xyoption {all}

\def \mv {{\mathfrak M_{v}}}  
\def \tv {{\Theta_{v}}}
\def \mw {{\mathfrak M_{w}}}
\def \tw {{\Theta_{w}}}

\newcommand{\comment}[1]{}
\newtheorem{theorem}{Theorem}
\newtheorem {lemma}{Lemma}
\newtheorem {question}{Question}

\theoremstyle{definition}

\theoremstyle {definition}
\newtheorem* {remark}{Remark}

\begin{document}
\baselineskip=16pt

\title[On Verlinde sheaves and strange duality]{On Verlinde sheaves and strange duality over elliptic Noether-Lefschetz divisors}
\author {Alina Marian}
\address {Department of Mathematics}
\address {Northeastern University}
\email{a.marian@neu.edu}
\author {Dragos Oprea}
\address {Department of Mathematics}
\address {University of California, San Diego}
\email {doprea@math.ucsd.edu}

\begin{abstract}
We extend results on generic strange duality for $K3$ surfaces by showing that the proposed isomorphism holds over an entire Noether-Lefschetz divisor in the moduli space of quasipolarized $K3$s. We interpret the statement globally as an isomorphism of sheaves over this divisor, and also describe the global construction over the space of polarized $K3$s. 
\end {abstract}
\maketitle

\addtocounter{section}{+1}
\subsection {Introduction} For a fixed polarized $K3$ surface $(X, H)$, let $v, w\in H^{\star}(X, \mathbb Z)$ be two primitive elements  which are orthogonal in the sense that $$\int_{X} v\cup w= 0.$$ Consider the moduli space $\mathfrak M_v$ of Gieseker $H$-stable sheaves $E$ on $X$ of Mukai vector $v$: $$\text{ch}(E)\sqrt{\text{Todd }(X)} =v.$$The Mukai vector $w$ induces a determinant line bundle $$\Theta_w\to \mathfrak M_v,$$ constructed in \cite {lp2}\cite {jli}. Specifically, if a universal family $\mathcal E\to \mv\times X$ is available, we set $$\Theta_w=\det \mathbf Rp_!(\mathcal E\otimes^{\mathbf L} q^{\star}F)^{-1},$$ for a complex $F\to X$ of Mukai vector $w$. Similarly we obtain the line bundle $\Theta_v\to \mw$. 

If $c_1(v\otimes w)\cdot H>0$, the defining equation of the locus $$\Theta=\{(E, F): \mathbb H^0(E\otimes^{\mathbf L} F)\neq 0\}\hookrightarrow \mv\times \mw$$ is a section of the line bundle $$\Theta_w\boxtimes \tv\to \mv\times \mw$$ inducing a map \begin{equation}\label{lpi}\mathsf D: H^0(\mv, \tw)^{\vee}\to H^0(\mw, \tv).\end{equation} According to Le Potier's strange duality conjecture \cite {lp1}, $\mathsf D$ is expected to be an isomorphism. 

In \cite{generic} we established the conjecture for generic polarized surfaces $(X, H)$ and for many pairs of Mukai vectors $(v, w)$ which satisfy $$c_1(v)=c_1(w)=H.$$ The proof involved degeneration to the locus of elliptic $K3$ surfaces with section and Picard rank $2$. 

Theorem \ref{p1} stated in Section $3$ strengthens this result by showing that the isomorphism holds over an entire Noether-Lefschetz divisor in the moduli space $\mathcal K_{\ell}$ of primitively {\it quasipolarized} $K3$ surfaces of degree $2\ell$. In Section $2$ we record basic properties of the Noether-Lefschetz divisor. In particular, while in \cite {generic} the elliptic fibrations are assumed to have only irreducible fibers, in the current setup, the fibers can be arbitrary.

In Section $4$, the duality is stated globally as an isomorphism of sheaves, the {\it Verlinde} sheaves, over the entire Noether-Lefschetz divisor. The Verlinde sheaves are also constructed more generally over the locus $\mathcal K_{\ell}^{\circ}\hookrightarrow \mathcal K_{\ell}$ of polarized $K3$s. It would be interesting to extend this construction to $\mathcal K_{\ell}$ in a suitable manner.
\addtocounter{section}{+1}
\subsection {The Noether-Lefschetz divisor $\mathcal P_1$}
Let  $(\mathcal X, \mathcal H)\to \mathcal K_{\ell}$ be the moduli stack of quasipolarized $K3$ surfaces $(X, H)$ of degree $H^2=2\ell$ with $\ell\neq 1$. 

We consider the Noether-Lefschetz loci of quasipolarized elliptically fibered $K3$ surfaces in ${\mathcal K}_{\ell}$. Specifically, for each $k >0$,  
we denote by $\mathcal P_k$ the Noether-Lefschetz stack parametrizing triples $(X, H, F)$ consisting of quasipolarized $K3$'s of degree $2\ell$, and divisor classes $F$ over $X$ satisfying
$$F^2 = 0, \, \, \, F \cdot H = k.$$ We claim that $$\mathcal P_1\hookrightarrow \mathcal K_{\ell}$$ is a substack of $\mathcal K_{\ell}$ parametrizing exactly the quasipolarized $K3$s which can be elliptically fibered with section, and with the quasipolarization a numerical section. This is expressed by the lemma below. The statement is standard, but a reference seemed difficult to find. 

\begin{lemma}\label{l1}
Let $(X, H)$ be a quasipolarized $K3$ surface of degree $2\ell$ with $\ell\neq 1$, and let $F$ be a divisor class on $X$ satisfying $$F^2 = 0, \, \, \, F \cdot H = 1.$$ Then
\begin{enumerate}
\item[(i)]
$F$ is effective and $\mathcal O (F)$ is globally generated;
\item[(ii)]
the induced map $\pi: X \to {\mathbb P}^1$ is an elliptic fibration with section $\sigma$, having $F$ as the fiber class;
\item [(iii)] the quasipolarization equals $H=\sigma + (\ell+1)F$; 
\item [(iv)] the class $F$ satisfying the two numerical assumptions above is unique.
\end{enumerate}
\end{lemma}

\proof Note first that $\chi ({\mathcal O} (F)) = 2.$ Since $-F \cdot H = -1,$ and $H$ is nef, $-F$ cannot be effective, so $$h^2 ({\mathcal O} (F)) = h^0 ({\mathcal O} (-F)) = 0, \, \, \, \text{and } \, \, \, h^0 ({\mathcal O} (F)) \geq \chi(\mathcal O(F))=2.$$ Thus $F$ is effective. 

We treat separately the two possibilities that $\mathcal O(F)$ be nef or not. First, if $\mathcal O(F)$ is nef, by the theorem of Piatetski-Shapiro and Shafarevich \cite {PS} there exists an elliptic fibration $$\pi:X\to \mathbb P^1$$ such that $F=mf$, where $f$ is the class of a fiber. In fact, $$F\cdot H=1\implies m=1,\,\,\,\,\, F=f, \,\,\,\,\,H\cdot f=1.$$ We next show that the fibration has a section. It is easy to check that the class $$\Sigma=H- (\ell+1) f$$ has self-intersection $-2$. Since $\chi(\mathcal O(\Sigma))=1,$ $\Sigma$ is either effective or anti-effective. In fact, $\Sigma$ is effective, since $\Sigma\cdot H>0$. Let $C$ be a curve in the linear series $\mathcal O(\Sigma)$. Now, for any component $R$ of a fiber we have $R\cdot f=0$ by Zariski's lemma, cf. III.8.2 \cite {barth}. Since $C\cdot f=1$, $C$ must have a component which intersects each fiber with multiplicity $1$. The other components of $C$ must be supported on  components of the fibers. The transversal component gives a section $\sigma$ of the elliptic fibration $\pi$. 

We now argue that $H=\sigma+(\ell+1)f$. From the above discussion, we already know that $$H=\sigma+mf+\sum m_i R_i$$ where $R_i$ are components of fibers and $m=\ell+1$. In fact, by absorbing other fiber classes into the constant $m$, we may assume $R_i$ are supported on fibers with two components or more. We have the following possibilities:
\begin {itemize}
\item [(i)] fibers of type $I_n$, consisting in a polygon of rational curves $C_1, \ldots, C_n$;
\item [(ii)] fibers of type $III$, consisting of $2$ rational curves $C_1, C_2$ meeting tangentially;
\item [(iii)] fibers of type $IV$ consisting of $3$ concurrent rational curves $C_1, C_2, C_3$;
\item [(iv)] fibers of type $I^{\star}_n$ which can be written as $$C_1+C_2+C_3+C_4+2(D_1+\ldots+D_n)$$ where $$C_1\cdot D_1=C_2\cdot D_1=C_3\cdot D_n=C_4\cdot D_n=1$$ and $D_i\cdot D_{i+1}=1$ for $1\leq i\leq n-1$;
\item [(v)] fibers of type $II^{\star}, III^{\star}, IV^{\star}$ corresponding to the graphs $E_6, E_7, E_8$. 
\end {itemize} Consider a fiber of type (i) and its contribution $\sum m_i C_i$ to the divisor $H$. We claim this contribution is a multiple of the fiber. Indeed, label the components so that $C_1$ intersects the section $\sigma$. Since $H\cdot C_i\geq 0$ for all $i$, we obtain the inequalities $$-2m_1+m_2+m_n\geq -1, \,-2m_2+m_1+m_3\geq 0, \,\ldots, -2m_n+m_1+m_{n-1}\geq 0.$$ If $-2m_1+m_2+m_n\geq 0$, then after adding the above inequalities, we conclude that we must have equality throughout. Thus $m_1=\ldots=m_n=m$ which shows that $\sum m_iC_i=mf$ as claimed. The case $$-2m_1+m_2+m_n=-1$$ is impossible. Indeed, since $$\sum_{k\neq 1} (-2m_k+m_{k-1}+m_{k+1})=-(-2m_1+m_2+m_n)=1$$ we conclude that for some index $k_0$ $$-2m_k+m_{k-1}+m_{k+1}=\begin{cases} 1\,\,\, \text{ if } k=k_0\\ 0\,\,\,\text { if } k\neq {1, k_0}.\end{cases}$$ This system is easily seen not to have any solutions. The remaining fiber types (ii)-(v) are entirely similar, and we will not verify them explicitly. In all cases, we find that $\sum m_i C_i$ must contribute a multiple of the fiber, hence $$H=\sigma + mf$$ for some integer $m$. In fact, $m=\ell+1$ by computing $H^2=2\ell$. This completes the proof when $\mathcal O(F)$ is nef. 

We assume now that $\mathcal O(F)$ is not nef and we will reach a contradiction. Then there exists an irreducible curve $\Gamma_1$ such that $$F\cdot \Gamma_1<0.$$ The curve $\Gamma_1$ is a component of an effective curve of class $F$ and furthermore $\Gamma_1^2<0$. Thus $\Gamma_1$ is a smooth rational curve on $X$. Let $H'$ be an ample class, and set $F_0=F$. The reflection of $F$ along $\Gamma_1$ then yields an effective class, cf. \cite {S}: $$F_1=F_0+(F_0\cdot \Gamma_1)\Gamma_1$$ which has the property that $$F_1^2=F_0^2=0,\,\,\, F_1\cdot H'<F_0\cdot H'.$$ If $F_1$ is not nef, then we continue the process reflecting along a smooth rational curve $\Gamma_2$. The process will eventually stop since $F_i\cdot H'$ is a decreasing sequence of non-negative integers. At the end, we find a nef line bundle $\mathcal O(F_k)$ of zero self-intersection, where $$F_k=F+(F_0\cdot \Gamma_1)\Gamma_1+(F_1\cdot \Gamma_2)\Gamma_2+\ldots + (F_{k-1}\cdot \Gamma_k)\Gamma_k.$$ Therefore $F_k=mf$, where $m\geq 0$ by nefness. In particular, $$F=mf+\sum n_i \Gamma_i$$ where $n_i=-F_{i-1}\cdot \Gamma_i>0.$ Using $F\cdot H=1$ we conclude $$m (H\cdot f)+\sum n_i (H\cdot \Gamma_i)=1.$$ Since $H$ is nef, the intersection numbers above are nonnegative. If $H\cdot f=0$, since $H^2>0$, by the Hodge index theorem we find $f^2\leq 0$. Since equality occurs, $f$ must be numerically trivial which is not the case since it intersects $H'$ nontrivially. Therefore $$H\cdot f=1,\,\, m=1,\,\,\,\,H\cdot \Gamma_i=0 \text{ for all } i.$$ The argument given in the nef case then shows that the elliptic fibration $\pi$ has a section $\sigma$, and $$H=\sigma+(\ell+1)f.$$ We conclude $$H\cdot \Gamma_i=\sigma\cdot \Gamma_i+(\ell+1) f\cdot \Gamma_i=0.$$ Thus either $\sigma\cdot \Gamma_i\leq 0$ or $f\cdot \Gamma_i\leq 0$. This means $\Gamma_i$ is contained in $\sigma$ or in the fiber $f$. The first case cannot occur since then $$\Gamma_i=\sigma\text{ and } \sigma\cdot \Gamma_i+(\ell+1) f\cdot \Gamma_i=0 \text{ shows } \ell=1$$ which is not allowed. Thus $\Gamma_i$ is a component of the fiber of $f$. However, in this case $f\cdot \Gamma_i=0$ by Zariski's lemma. Since $$F=f+\sum n_i\Gamma_i$$ has zero self intersection, we find $$(\sum n_i \Gamma_i)^2=0,$$ where $\Gamma_i$ are components of the fiber. This yields $\sum n_i \Gamma_i=nf$ for some integer $n$, again by Zariski's lemma. Thus $F=(n+1)f$, and since $F\cdot H=1$ then $F$ is the fiber class. 

Finally, we establish the uniqueness of $F$ as claimed in (iv). If $F'$ is another class with $$F'^2=0, \,\, F'\cdot H=1$$ then we can write $$F'=a\sigma+R$$ where $R$ is supported on components of fibers. We have $R\cdot f=0$ and $$F'\cdot H=\left(a\sigma+R\right)\cdot (\sigma+(\ell+1)f)=1\implies R\cdot \sigma=1-a(\ell-1).$$ In addition $$F'^2=0\implies -2a^2+2a (R\cdot \sigma)+R^2=0.$$ This yields $$R^2=-2a+2a^2(\ell+1).$$By Zariski's lemma, $R^2\leq 0$, which implies $a=0$. Furthermore, we obtain $R^2=0$, showing that $R=mf$, again by Zariski's lemma. Moreover, $R\cdot \sigma = 1$ hence $m=1$. Therefore $F'=f$, proving uniqueness.
\qed
\addtocounter{section}{+1}
\subsection {Strange duality along $\mathcal P_1$} 
For $(X, H)\in \mathcal P_1$, we consider the orthogonal Mukai vectors  \begin{equation}\label{mvec}v=r+H+a\,[\text{pt}], \,\,\, w=s+H+b\,[\text{pt}]\end{equation} with $r, s\geq 3$, satisfying further 
\begin{equation}\label{ii}\langle v, v\rangle + \langle w, w\rangle \geq 2(r+s)^2.\end{equation}
We form the moduli spaces of sheaves\footnote {Stability of the sheaves in $\mv$ and $\mw$ is with respect to a polarization which is suitable in the sense of Friedman \cite {F}. As shown in the appendix of \cite {generic}, this choice of polarization is in fact irrelevant under the stronger assumptions that $$\langle v, v\rangle \geq 2(r-1)(r^2+1),\,\,\, \langle w, w\rangle \geq 2(s-1)(s^2+1).$$ Indeed, in this case, the different moduli spaces are birational in codimension $1$.}  $\mv$ and $\mw$ together with the corresponding theta line bundles. 

Under these conditions, in \cite{generic}, the strange duality map $$\mathsf D:\, H^0(\mv, \tw)^{\vee}\to H^0(\mw, \tv)$$ was proven to be an isomorphism over the open sublocus of $\mathcal P_1$ consisting of surfaces with Picard rank $2$. We now show:
\begin {theorem} \label{main}For the Mukai vectors $v, w$ specified above, the strange duality map $\mathsf D$ is an isomorphism over the entire Noether-Lefschetz divisor $\mathcal P_1$. 
\label{p1}
\end {theorem}

\proof For surfaces in $\mathcal P_1$ of Picard rank larger than 2, the elliptic fibration has finitely many reducible fibers. Fourier-Mukai functors were studied in this setting in \cite{hms2}. Specifically, let $$\pi: X \to {\mathbb P}^1$$ be any quasipolarized elliptically fibered $K3$ surface with section class $\sigma$ and fiber class $f$. Consider the product $Y=X \times_{{\mathbb P}^1} X$ with projections $p$ and $q$ to the two factors, and let $$\Delta \subset X \times_{{\mathbb P}^1} X$$ be the diagonal. The $\pi$-relative Fourier-Mukai functor 
$${\mathsf S}: {\mathbf D} (X) \longrightarrow {\mathbf D}(X)$$
with kernel 
$${\mathcal P} = {\mathcal I}_{\Delta}\otimes {\mathcal O} (p^{\star} \sigma + q^{\star} \sigma)$$
is an equivalence of bounded derived categories of coherent sheaves by Proposition 2.16 of \cite{hms2}. As $(X, H)$ is in $\mathcal P_1$, by Lemma \ref{l1} $$c_1(v)=c_1(w) =\sigma+(\ell+1)f.$$ Along the lines of \cite{bridgeland}, we shall prove shortly that the Fourier-Mukai transform ${\mathsf S}$ induces a birational morphism, regular in codimension $1$, between the moduli spaces ${\mathfrak M}_v$ and ${\mathfrak M}_w$ on the one hand,  and the Hilbert schemes of $d_v$ respectively $d_w$ points on $X$ on the other: 
$$\Psi_v: \mathfrak M_v\dashrightarrow X^{[d_v]}, \,\,\,\,\,\ \Psi_w: \mathfrak M_w\dashrightarrow X^{[d_w]}.$$

Assuming this for the moment, we explain how to complete the proof of Theorem \ref{main}, much as in \cite {generic}. We determine first the exact numerics of the transformation $\mathsf S$ by a cohomological Fourier-Mukai calculation. Let $V \in {\mathbf D} (X)$ be any complex of rank $r$, Euler characteristic $\chi$, and first Chern class $$c_1(V) = k \sigma + m f,$$ for integers $k$ and $m$. Recalling $p$ and $q$ are the projections from $Y=X\times_{\mathbb P^1} X$, we have \begin{eqnarray*}\det \mathsf S(V)&=&\det \mathbf Rq_{\star}(\mathcal P\otimes p^{\star}V)=\det \mathbf R q_{\star}(I_{\Delta}\otimes p^{\star}V(\sigma)\otimes q^{\star} \mathcal O(\sigma))\\ &=& \det \mathbf Rq_{\star}(I_{\Delta}\otimes p^{\star} V(\sigma))\otimes \mathcal O(\sigma)^{\chi(V|_{f})}\\ &=& \det \mathbf Rq_{\star}(p^{\star} V(\sigma))\otimes \det\mathbf Rq_{\star}(\mathcal O_{\Delta}\otimes p^{\star}V(\sigma))^{-1}\otimes \mathcal O(k\sigma)\\&=& \det \mathbf Rq_{\star}(p^{\star} V(\sigma))\otimes \det V(\sigma)^{-1}\otimes \mathcal O(k\sigma)\\&=&\det \mathbf Rq_{\star}(p^{\star} V(\sigma))\otimes \mathcal O(-r\sigma - mf).\end{eqnarray*}
To calculate the first term, it is more convenient to work on the product $$j:Y\hookrightarrow X\times X.$$ Let $\bar p, \bar q$ denote the two projections from $X\times X$, and let ${\text{pr}} = \pi \times \pi : \, X \times X \to {\mathbb P}^1 \times {\mathbb P}^1.$ Observing that $$j_{\star}\mathcal O_Y=\text{pr}^{\star} \mathcal O_{\Delta/\mathbb P^1\times \mathbb P^1}=\text{pr}^{\star}(\mathcal O_{\mathbb P^1\times \mathbb P^1}-\mathcal O_{\mathbb P^1}(-1)\boxtimes \mathcal O_{\mathbb P^1}(-1))=\mathcal O_{X\times X}-\bar p^{\star}\mathcal O(-f)\otimes \bar q^{\star} \mathcal O(-f),$$ we calculate \begin{eqnarray*} \det \mathbf Rq_{\star}(p^{\star} V(\sigma))&=&\det \mathbf R{\bar q}_{\star}\,\left ({\bar p}^{\star} V(\sigma)\otimes j_{\star}\mathcal O_Y\right )\\ &=& \det \mathbf R{\bar q}_{\star}\,(\bar p^{\star} V(\sigma))\otimes \det \mathbf R\bar q_{\star}\,\left (\bar p^{\star} V(\sigma)\otimes \bar p^{\star}\mathcal O(-f)\otimes \bar q^{\star}\mathcal O(-f)\right )^{-1}\\&=&\det (\mathbf R\bar q_{\star}\, \left (\bar p^{\star} V(\sigma-f))\otimes \mathcal O(-f)\right )^{-1}\\&=&\mathcal O(-f)^{-\chi(V(\sigma-f))}=\mathcal O((\chi-2r+m-3k)f).\end{eqnarray*} To summarize, we obtained \begin{eqnarray}
 \label{c1}
 \det {\mathsf S} (V) = \mathcal O(-r \,\sigma + (\chi - 2r -3k) \,f). \nonumber
 \end{eqnarray}

\vskip.1in

Now let $E$ and $F$ be stable sheaves whose Mukai vectors $v$ and $w$ are given by \eqref{mvec}. By the preceding calculation 
$$\det {\mathsf S} (E^{\vee}) = \mathcal O(-r\sigma + (a-r+3) f),$$
$$\det {\mathsf S} (F) = \mathcal O( -s \sigma + (b-s -3) f). $$ Assuming the birational isomorphism with the Hilbert scheme, for generic $E$ and $F$ we therefore have that 
\begin{equation}
\label{fm1}
{\mathsf S} (E^{\vee}) = I_Z \otimes {\mathcal O} (r\sigma - (a-r+3) f) [-1], 
\end{equation}
\begin{equation}
\label{fm2}
{\mathsf S} (F) = I_{W}^{\vee}  \otimes {\mathcal O} (-s \sigma + (b-s -3) f),
\end{equation}
where $Z$ and $W$ are zero dimensional subschemes of lengths $d_v$ and $d_w$ respectively. In fact, we will only explain the first equality below; the second can be deduced from the first by Grothendieck duality as in Proposition $2$ of \cite {generic}. 

We finally calculate
\begin{eqnarray*}
\mathbb H^0 (E \otimes^{\mathbf L} F) &=& \text{Hom}_{{\mathbf D}(X)} (E^{\vee},\, F) =  \text{Hom}_{{\mathbf D}(X)} \left ({\mathsf S} (E^{\vee}), \,{\mathsf S} (F) \right ) \\
& = & \text{Ext}^1 (I_Z \otimes L, \, I_W^{\vee} )= \text{Ext}^1 (I_W^{\vee}, \, I_Z \otimes L)^{\vee} \\
& = & \mathbb H^1 (I_W \otimes^{\mathbf L} I_Z \otimes L)^{\vee}.
\end{eqnarray*}
On the third line, using \eqref{fm1} and \eqref{fm2}, we have set 
$$L = {\mathcal O} \left (( r+s) \sigma +(r+s - a- b) f \right).$$
The orthogonality condition $$H^2 = - r b -s a$$ for the Mukai vectors $v$ and $w$ together with the bound \eqref{ii} on the dimensions $d_v$ and $d_w$ ensure that
$ -a - b > r +s,$ so the line bundle $L$ is big and nef, without higher cohomology on $X$.  

\vskip.1in

Thus, under the birational map
$$\Psi_v \times \Psi_w: \, \, {\mathfrak M}_v \times {\mathfrak M}_w \dasharrow X^{[d_v]} \times X^{[d_w]}$$ the two theta divisors $$\Theta = \{(E, \, F): \, \, \mathbb H^0 (E \otimes^{\mathbf L} F) \neq 0 \} \subset {\mathfrak M}_v \times {\mathfrak M}_w,$$ and 
$$\theta_L = \{(I_Z, \, I_W): \, \, \mathbb H^0 (I_Z \otimes^{\mathbf L} I_W \otimes L) \neq 0 \}  \subset X^{[d_v]} \times X^{[d_w]}$$ coincide. The line bundles $\Theta_w, \, \Theta_v$ on the two higher-rank moduli spaces and $L^{[d_v]}, \, L^{[d_w]}$ on the two Hilbert schemes are also identified. As explained in \cite{generic}, for line bundles $L$ without higher cohomology, $\theta_L$ is known to induce an isomorphism \begin{equation}\label{sds}H^0 (X^{[d_v]}, \, L^{[d_v]} )^{\vee} \longrightarrow H^0 ( X^{[d_w]}, \, L^{[d_w]} ).\end{equation} Therefore, under the identifications above, $\Theta$ also induces the isomorphism of equation \eqref{lpi}:
$$\mathsf D: H^0 ( {\mathfrak M}_v, \, \Theta_w)^{\vee} \longrightarrow H^0 ({\mathfrak M}_w, \, \Theta_v).$$

\vskip.1in

We turn now to the proof that $\Psi_v$ is an isomorphism in codimension 1, which was given for a surface $\pi: X \to {\mathbb P}^1$ with irreducible fibers in \cite{bh}, \cite{generic}. We thus take up the case when the fibration has at least one reducible fiber. We shall explain that the {\it inverse} $$\Psi_v^{-1}: X^{[d_v]} \dasharrow {\mathfrak M}_v$$ is a regular embedding defined on a subscheme $U \subset X^{[d_v]}$ with $\text{codim} \, (X^{[d_v]} \setminus U) \geq 2.$ The same is then true about $\Psi_v$ on ${\mathfrak M}_v.$ Indeed, if this were not the case, as the two moduli spaces are holomorphic symplectic, $\Psi_v$ would at least admit by \cite{H} an extension $\overline{\Psi}_v$ to a regular embedding defined away from codimension 2 on ${\mathfrak M}_v$. Thus $\overline{\Psi}_v$ would extend over a divisorial locus $D \subset {\mathfrak M}_v$ where the original map $\Psi_v$ is assumed undefined. But then $${\overline{\Psi}}_v (D) \subset X^{[d_v]} \setminus U,$$ a contradiction as the latter has codimension 2 in $X^{[d_v]}.$

\comment{It is simpler to explain that the inverse 
$$\Psi_v^{-1}: X^{[d_v]} \dasharrow {\mathfrak M}_v$$ is defined away from codimension 2 on the Hilbert scheme when $r\geq 3$. As both $\mathfrak M_v$ and $X^{[d_v]}$ are holomorphic symplectic, the same is then true about $\Psi_v$. Indeed, $\Psi_v$ admits an extension to a regular map $f$ over a maximal open set $U\hookrightarrow \mathfrak M_v$ whose complement has codimension at least $2$. Furthermore, $f$ is an embedding over $U$, cf. \cite {H}. Now assuming the Fourier-Mukai $\Psi_v$ is undefined along $D\subset \mathfrak M_v$, then $$\mathfrak M_v\setminus D\subset U.$$ The inverse $\Psi_v^{-1}$ is defined along $\Psi_v(\mathfrak M_v\setminus D).$ The complement of this set in the Hilbert scheme has codimension $2$. However $$f(D\cap U)\subset X^{[d_v]}\setminus \Psi_v(\mathfrak M_v\setminus D)$$ and since $f$ is an embedding, $D$ must have codimension at least $2$ as well. }

We are thus left to analyze the domain of $\Psi_v^{-1}$. The inverse is a Fourier-Mukai transform whose kernel is a complex $\mathcal Q[1]$ over $X\times_{\mathbb P^1} X$. We write $\mathsf T$ for the Fourier-Mukai transform with kernel $\mathcal Q$ so that $$\mathsf S\circ \mathsf T=[-1],\,\, \mathsf T\circ \mathsf S=[-1].$$ We claim that for generic $Z$, the sheaf  $$M=I_Z \otimes {\mathcal O} (r \sigma - (a-r+3) f)$$ is $WIT_0$ for the kernel $\mathcal Q$. Its transform is then a stable torsion free sheaf in $\mv$, cf. Section $7$ of \cite {bridgeland}. To prove the claim, we adapt arguments of \cite {bridgeland}, as follows. On general grounds, cf. Lemma $6.1$ in \cite {bridgeland}, there is a short exact sequence $$0\to A\to M\to B\to 0$$ where $A$ is $\mathsf T$-$WIT_0$, while $B$ is $\mathsf T$-$WIT_1$. We prove that $B=0$, following Lemma $6.4$ in \cite {bridgeland}. Assuming otherwise, we have $\mathsf T(B)\neq 0$, and therefore there exists $x\in X$ and a non-zero morphism $$\mathsf T^1(B)\to \mathbb C_{x}.$$ Note however that $$\mathbb C_x=\mathsf T^1(I_{x}(o)),$$ where $I_x$ is the ideal sheaf of the point $x$ in its fiber, and $o$ denotes the intersection of the fiber through $x$ with the section.  In fact, $I_x(o)=\mathsf S^0(\mathbb C_x)$, by Lemma $6.3.7$ of \cite {C}. By Parseval, we now obtain a non-zero morphism $$M\to B\to I_x(o).$$ This morphism must factor through the restriction of $M$ to the fiber $C$ through $x$, yielding a non-zero map $$I_Z|_{C}\otimes \mathcal O(ro)\to I_x(o).$$ Thus it suffices to show $$\text{Hom}_C(I_Z|_{C}\otimes \mathcal O((r-1)o), I_x)=0.$$ We prove this is the case for $r\geq 3$ and subschemes $Z$ such that 
\begin {itemize}
\item [(i)] $Z$ intersects any smooth fiber in at most two points; 
\item [(ii)] $Z$ intersects any singular fiber in at most one point which is not a node or a cusp (if the fiber is irreducible) or does not lie on at least two irreducible components. 
\end {itemize}
This locus has complement of codimension $2$ in the Hilbert scheme of $X$. 

When $C$ is a smooth fiber, $\zeta=Z\cap C$ has length at most equal to $2$, by (i). Then $$I_Z|_{C}=I_{\zeta/C}\oplus T$$ where $T$ is a torsion sheaf supported at $\zeta$. This can be seen by restricting the ideal sequence of $Z$ to the curve $C$. In fact, the same statement also holds when $C$ is singular, as $Z$ is subject to (ii). When $C$ is smooth, it suffices therefore to prove $$\text{Hom}_C(I_{\zeta/C}((r-1)o), I_x)=0\iff H^0(\mathcal O_C(-(r-1)o+\zeta-x))=0.$$ Since for $r\geq 3$ the degree is negative, the conclusion follows. When $C$ is a singular fiber, the scheme $\zeta=Z\cap C$ has length $1$. We show $$\text{Hom}_C(I_{\zeta/C}((r-1)o), \mathcal O_C)=0\text{ which gives } \text{Hom}_C(I_{\zeta/C}((r-1)o), I_x)=0.$$ Indeed, by duality, this is the same as proving $$H^1(I_{\zeta/C}((r-1)o)))=0.$$
Here we used that the dualizing sheaf of $C$ is trivial. Assume first $\zeta\neq o$. From the exact sequence $$0\to I_{\zeta/C}(o)\to I_{\zeta/C}((r-1)o)\to \mathbb C_{o}^{r-2}\to 0$$ we see it suffices to show $$H^1(I_{\zeta/C}(o))=0.$$ Next, from the exact sequence $$0\to \mathcal O(-o)\to \mathcal O\to \mathbb C_o\to 0$$ we conclude $$H^0(\mathcal O(-o))=0, \,\, H^1(\mathcal O(-o))=\mathbb C\implies H^0(\mathcal O(o))=\mathbb C,\,\, H^1(\mathcal O(o))=0.$$ The exact sequence $$0\to I_{\zeta/C}(o)\to \mathcal O_C(o)\to \mathbb C_{\zeta}\to 0$$ and the fact that $$H^0(\mathcal O_C(o))\to \mathbb C_{\zeta}$$ is an isomorphism for $\zeta\neq o$ yield $H^1(I_{\zeta/C}(o))=0,$ as claimed. The vanishing of higher cohomology also holds for $\zeta=o$ since $H^1(\mathcal O((r-2)o))=0$. This completes the proof. \qed

\addtocounter{section}{+1}
\subsection{The Verlinde sheaves} We will reinterpret Theorem \ref{main} as giving an isomorphism of sheaves defined over the divisor $\mathcal P_1$ in the moduli space of quasipolarized $K3$s. 

\subsubsection {Construction} For a fixed integer $n$, we may consider over  ${\mathcal K}_{\ell}$ the relative Hilbert scheme of $n$ points $$\pi: {\mathcal X}^{[n]} \to {\mathcal K}_{\ell}, $$ viewed as the relative moduli stack of rank 1 torsion free sheaves of trivial determinant and second Chern number $-n$. 

More generally, to consider spaces of higher rank sheaves as the $K3$ surface varies in moduli, we restrict attention to the open substack ${\mathcal K}_{\ell}^{\circ} \subset {\mathcal K}_{\ell}$ where the line bundle $\mathcal H$ over the universal surface $$\pi:\mathcal X\to \mathcal K_{\ell}$$ is ample. We construct $$M[v]\to \mathcal K_{\ell}^{\circ},$$ the moduli space of $\mathcal H$-semistable sheaves 
with rank $r$, determinant $d {\mathcal H}$ and Euler characteristic $a-r$ over the fibers of $\pi:  {\mathcal X}^\circ\to {\mathcal K}_{\ell}^\circ.$ 

The construction of the theta bundles over $M[v]$ is subtler. To start, let $$\pi: \mathcal X_1^{\circ} \to {\mathcal K}^{\circ}_{\ell, 1}$$ be the universal family over the moduli stack ${\mathcal K}^{\circ}_{\ell, 1}$ of polarized $K3$s with a marked point.  It has a canonical section $$\sigma: {\mathcal K}^{\circ}_{\ell, 1} \to {\mathcal X}_1^{\circ}.$$  
\comment{After contracting the $(-2)$-curves, this yields a section $$\tilde \sigma:\mathcal K_{\ell, 1}\to \widetilde {\mathcal X}_1.$$ }Let 
\begin{eqnarray}
\nonumber
\label{basicvecglobal}
\mathcal V&= & (r-d)\, {\mathcal O} + d\, {\mathcal H} + \alpha\, {\mathcal O}_ { \sigma}, \\ \mathcal W&=& (s-e)\, {\mathcal O} + e\, {\mathcal H} + \beta \, {\mathcal O}_{ \sigma}, \nonumber
\end{eqnarray}
be classes in the $K$-theory of $ {\mathcal X}_1^{\circ}.$ Over a fixed marked polarized $K3$ surface $(X, H, p)$, they have the Mukai vectors $$v=r+dH+a[\text{pt}],\,\,\, w=s+eH+b[\text{pt}],$$ for 
$$\alpha= a - r - \frac{d H^2}{2},$$
$$ \beta= b - s - \frac{e H^2}{2}.$$ We further denote as $$\pi_v: M[v]_1 \longrightarrow {\mathcal K}^{\circ}_{\ell, 1}$$ the relative moduli space of stable sheaves of type $v$ over the fibers of $ \pi:  {\mathcal X}_1^{\circ}\to \mathcal K_{\ell, 1}^{\circ}.$ The class $\mathcal W$ induces standardly a determinant line bundle $$\overline{\Theta}_w \to M[v]_1,$$ via descent from $$\mathcal Q \to M[v]_1,$$ where $\mathcal Q$ is an open subscheme of a suitable quot scheme. Explicitly, over $\mathcal Q$, we have $$\overline{\Theta}_w=\det \mathbf Rp_{!}(\mathcal E\otimes q^{\star} \mathcal W)^{-1}$$ for the universal quotient sheaf  $\mathcal E \to \mathcal Q\times_{\mathcal K^{\circ}_{\ell, 1}}{\mathcal X}_1^{\circ}.$ The fiber of the forgetful map $$M[v]_1 \, \to \, M[v]$$ over a point $(X, \, H, \, E \to X) \in M[v]$ is the surface $ X$. To describe the restriction of $\overline{\Theta}_w$ to this fiber, we let $\Delta \subset  X \times X $ be the diagonal and denote by $p, q$ the projections from $X \times  X$ to the two factors.
Then $$\left . \overline{\Theta}_w \right |_{ X} = {\det} \, {\mathbf R}p_{\star} \left ( q^{\star} E \otimes ( (s-e) {\mathcal O} \oplus q^{\star} (e H) \oplus \beta\, {\mathcal O}_{\Delta} ) \right )^{-1} = \det E^{-\beta} =  H^{-\beta d}. $$ 
We conclude that the product line bundle 
\begin{equation}
\label{norm}
\overline{\Theta}_w \otimes \pi_{v}^{\star} {{\mathcal H}}^{\,\beta d} \, \, \, \text{on} \, \, \, M[v]_1
\end{equation}
restricts trivially to the fibers of the map $$M[v]_1 \longrightarrow  M[v]$$ forgetting the marking. By the seesaw lemma, the product \eqref{norm} is in fact the pullback to $M[v]_1$ of a line bundle ${\Theta}_w\to M[v]:$ $$\overline{\Theta}_w \otimes \pi_{v}^{\star} {{\mathcal H}}^{\,\beta d}=\text{pr}^{\star}{ \Theta}_w.$$
While the determinant line bundle $\Theta_w$ is uniquely defined for a fixed $K3$ surface, over the relative moduli space $M[v]$, $\Theta_w$ depends on choice of $\mathcal H$, and therefore can be canonically defined only up to tensoring by line bundles pulled back from ${\mathcal K}_{\ell}^{\circ}$. 

\begin{remark} The same construction gives the theta line bundle on the relative moduli space ${\mathcal S}{\mathcal U}_g (r) \longrightarrow M_g$ of semistable rank $r$ bundles with trivial determinant over smooth curves of genus $g$. They are naturally defined on the basechanged moduli space
 $${\mathcal S}{\mathcal U}_{g,1} (r)=  {\mathcal S}{\mathcal U}_g (r) \times_{M_g} M_{g,1} \longrightarrow M_{g,1},$$ relative to the $K$-theory class $${\mathcal O} + (g-1) {\mathcal O}_{\sigma} $$ on the universal curve ${\mathcal C} \to M_{g,1},$ and are then seen to be pulled back under the forgetful map $${\mathcal S}{\mathcal U}_{g,1} (r) \to  {\mathcal S}{\mathcal U}_g (r).$$ Pushing forward the $k$-tensor powers of the theta line bundles to $M_g$, we obtain the Verlinde bundles $$\mathcal V_{r,k}\to M_{g}.$$ Their first Chern classes remain unknown in general.
\end{remark} 
\vskip.1in

\subsubsection {Global strange duality.} Over ${\mathcal K}_{\ell}^{\circ} $ we define now the Verlinde complexes
\begin{equation}
\label{versheaves}
{\mathbf W}= {\mathbf R}{\pi_v}_{\star} \Theta_w, \, \, \, \, \, {\mathbf V} = {\mathbf R}{\pi_w}_{\star} \Theta_v.
\end{equation} Consider the fiber product $$\pi: M[v] \times_{{\mathcal K}^{\circ}_{\ell}} M[w] \to {\mathcal K}^{\circ}_{\ell},$$ endowed with the canonical Brill-Noether locus,
\begin{equation} 
\label{thetadef}
\Theta = \{ (X, H, E, \, F) \, \, \, \text{so that} \, \, \, \mathbb H^0 (X, E \otimes^{\mathbf L} F) \neq 0\} \subset M[v] \times_{{\mathcal K}^{\circ}_{\ell}} M[w].
\end{equation}
One expects $\Theta$ to be a divisor. This was established in \cite{generic} when $v$ and $w$ satisfy $$c_1 (v) = c_1 (w) = {\mathcal H}.$$ The corresponding line bundle, also denoted for simplicity as $\Theta$, is in any case always defined on the product space, and splits by the seesaw lemma as
\begin{equation}
\Theta  \simeq \Theta_w \boxtimes \Theta_v.
\end{equation}
The above equation is correct up to a twist $\mathcal T\to \mathcal K_{\ell}^{\circ}$ which will be found explicitly below, and which for now we absorb into any one of the theta bundles. The two line bundles $\Theta_w$ and $\Theta_v$ are ambiguous up to reverse twistings by a line bundle from ${\mathcal K}_{\ell}^{\circ}$,
$$\left ( \Theta_v, \, \Theta_w\right ) \,  \sim \, \left ( \Theta_v \otimes \pi_w^{\star} {\mathcal L}, \, \Theta_w \otimes \pi_v^{\star} {\mathcal L^{-1}} \right ), \, \, \, \text{for} \, \, \, {\mathcal L} \in \, \text{Pic} \, {\mathcal K}_{\ell}^{\circ},$$ while $\Theta$ is canonical. Pushing forward the canonical theta line bundle via 
$\pi,$ we get
\begin{equation}
\label{pushtheta}
{\mathbf R}\mathbf {\pi}_{\star} \Theta \simeq  {\mathbf W} \otimes^{\mathbf L} {\mathbf V},
\end{equation}
and the above ambiguity carries over to the Verlinde complexes ${\mathbf W}$ and ${\mathbf V}.$
The divisor \eqref{thetadef} then induces a morphism $$\mathsf D:\mathbf W^{\vee}\to \mathbf V.$$ 
In \cite {generic}, also having assumed that $$\chi(v), \chi(w)\leq 0,$$ we showed that  over a Zariski open subset of $\mathcal K_{\ell}^{\circ}$, the higher cohomology sheaves vanish while $\mathcal H^{0}(\mathsf D)$ induces an isomorphism between the zeroth cohomology sheaves.

\begin {remark} Even though not necessary for our argument, let us determine the twist $\mathcal T\to \mathcal K_{\ell}^{\circ}$ in the decomposition \begin{equation}\label{pullback}\Theta=\Theta_w\boxtimes \Theta_v\otimes \text{pr}^{\star} \mathcal T\end{equation} over $M[v]\times_{\mathcal K_{\ell}^{\circ}} M[w]$. Above, we absorbed this twist into the Verlinde complexes, for the ease of exposition. 

First, we may pass to the moduli stack $\mathcal M[v]$ and $\mathcal M[w]$ of all sheaves over $X$, without changing the above equations.  We let $${\mathcal E} \to \mathcal M[v]_1 \times_{{\mathcal K}_{\ell, 1}^{\circ}} {\mathcal X}_1^{\circ}, \, \, \, \, \, \,    {\mathcal F} \to \mathcal M[w]_1 \times_{{\mathcal K}_{\ell, 1}^{\circ}} {\mathcal X}_1^{\circ}$$ be the universal families of sheaves,
and further set, on the same product spaces, 
$$\overline{\mathcal E} = \mathcal E - \text{pr}_2^{\star} \mathcal V, \, \,  \, \,  \, \, \, \, \overline{\mathcal F} = \mathcal F -\text{pr}^{\star}_2 \mathcal W.$$ Considering now the triple product 
$$\mathcal M[v]_1 \times_{{\mathcal K}_{\ell, 1}^{\circ}} \mathcal M [w]_1 \times_{{\mathcal K}_{\ell, 1}^{\circ} } {\mathcal X}_1^{\circ},$$ we calculate $$\Theta \otimes \Theta_v^{-1} \otimes \Theta_w^{-1}$$ as the pushforward
$$\left({\det} {\mathbf Rp_{12}}_{\star} \left ( p_{13}^{\star} \mathcal E \otimes^{\mathbf L} p_{23}^{\star} \mathcal F -  p_{13}^{\star} \mathcal E \otimes^{\mathbf L} p_3^{\star} {\mathcal W} -  p_{23}^{\star} \mathcal F \otimes^{\mathbf L} p_3^{\star} {\mathcal V}  \right )\right)^{-1} \otimes \text{pr}^{\star}{\mathcal H}^{-d\beta - e\alpha}$$ $$= \left({\det} {\mathbf Rp_{12}}_{\star} \left ( p_{13}^{\star} \overline{\mathcal E} \otimes^{\mathbf L} p_{23}^{\star} \overline{\mathcal F} - p_3^{\star} (\mathcal V \otimes^{\mathbf L} \mathcal W)\right)\right)^{-1} \otimes \text{pr}^{\star}{\mathcal H}^{-d \beta  - e \alpha },$$
where $\mathcal H \to \mathcal K_{\ell, 1}^{\circ}$ is viewed on $\mathcal M[v]_1 \times_{{\mathcal K}_{\ell, 1}^{\circ}} \mathcal M [w]_1$ via pullback by the natural projection  $$\text{pr}:\mathcal M[v]_1\times_{\mathcal K_{\ell,1}^{\circ}}\mathcal M[w]_1\to\mathcal K_{\ell,1}^{\circ}.$$ We apply Grothendieck-Riemann-Roch to compute $$\text{ch} \, {\mathbf Rp_{12}}_{\star} \left ( p_{13}^{\star} \overline{\mathcal E} \otimes^{\mathbf L} p_{23}^{\star} \overline{\mathcal F} \right ).$$ By construction, $\text{ch}\,  \overline{\mathcal E}$ and $\text{ch}\,  \overline {\mathcal F}$ restrict trivially over the fibers of 
$$p_{12}:\mathcal M[v]_1 \times_{{\mathcal K}_{\ell, 1}^{\circ}} \mathcal M [w]_1 \times_{{\mathcal K}_{\ell, 1}^{\circ} } {\mathcal X}_1^{\circ}\to \mathcal M[v]_1\times_{\mathcal K_{\ell,1}^{\circ}}\mathcal M[w]_1.$$ The Chern character of the pushforward above is thus supported in codimension $2$ or higher, and therefore gives $$ {\det} \,{ \mathbf Rp_{12}}_{\star} \left ( p_{13}^{\star} \overline{\mathcal E} \otimes^{\mathbf L} p_{23}^{\star} \overline{\mathcal F} \right ) = \mathcal O.$$
Recalling the morphism $\pi: {\mathcal X}_1^{\circ} \to \mathcal K_{\ell, 1}^{\circ}$ which describes the universal surface,  we find that  
\begin{eqnarray*}
\Theta &\otimes& \Theta_v^{-1} \otimes \Theta_w^{-1}  =  \det  {\mathbf Rp_{12}}_{\star} \left [p_3^{\star} (\mathcal V \otimes^{\mathbf L} \mathcal W) \right ]  \otimes  \text{pr}^{\star} {\mathcal H}^{-d \beta  - e \alpha } \\
& = & \text{pr}^{\star} \left(\det \mathbf R \pi_{\star} \left ( \mathcal V \otimes^{\mathbf L} \mathcal W \right ) \otimes  {\mathcal H}^{-d \beta  - e \alpha } \right) \\
& = & \text{pr}^{\star} \left(\det \mathbf R \pi_{\star} \left [\left ( (r-d){\mathcal O} + d{\mathcal H} +\alpha{\mathcal O}_ { \sigma}\right )  \otimes^{\mathbf L} \left ((s-e){\mathcal O} + e{\mathcal H} + \beta{\mathcal O}_{ \sigma}\right )  \right ] \otimes  {\mathcal H}^{-d \beta  - e \alpha }\right) \\
& = &  \text{pr}^{\star}\left (  \lambda^{-(r-d)(s-e)} \otimes \left ( \det \pi_{\star} {\mathcal H} \right )^{e(r-d) + d(s-e)} \otimes \left ( \det \pi_{\star} \mathcal H^2 \right )^{de} \right).
\end{eqnarray*}
Here, we wrote $$\lambda=(\det {\mathbf R}\pi_{\star}\mathcal O_{\mathcal X})^{-1}\to \mathcal K_\ell$$ for the Hodge bundle.
This yields 
$$\mathcal T =  \lambda^{-(r-d)(s-e)} \otimes \left ( \det \pi_{\star} {\mathcal H} \right )^{e(r-d)+d(s-e)} \otimes \left ( \det \pi_{\star} \mathcal H^2 \right )^{de}.$$
\end {remark}
\subsubsection {Extensions of the Verlinde sheaves and desiderata}
We now turn our attention to the locus of elliptic $K3$ with section, where the Verlinde sheaves  and the isomorphism $\mathsf D$ can be extended from $$\mathcal P_1^{\circ}=\mathcal P_1\cap \mathcal K_{\ell}^{\circ}$$ to all of $\mathcal P_1$ by the results of Section $3$, as we now explain. 

The universal data over $\mathcal P_1$ consists of the triple $$(\mathcal X, \mathcal H, \mathcal F)\to \mathcal P_1,$$ where $\mathcal F$ denotes the universal fiber class of the elliptic fibration. We consider the line bundle $$\mathcal L=\mathcal H^{r+s}\otimes \mathcal O(\mathcal F)^{-(r+s)\,\ell-a-b},$$ which restricts over each $(X, H, F)$ to $$L=\mathcal O((r+s)\sigma+(r+s-a-b)f).$$ In the product of Hilbert schemes we have the universal theta divisor $$\theta=\{(X, Z, W): \mathbb H^0(X, \,\,I_Z\otimes^{\mathbf L} I_W\otimes \mathcal L|_{X})\neq 0\}\subset \mathcal X^{[d_v]}\times_{\mathcal P_1}\mathcal X^{[d_w]}.$$ To write the corresponding line bundle, we denote by $${\mathcal Z} \subset {\mathcal X}^{[d_v]} \times_{{\mathcal K}_{\ell}} {\mathcal X}, \, \, \, \, \, \, {\mathcal W} \subset {\mathcal X}^{[d_w]} \times_{{\mathcal K}_{\ell}} {\mathcal X},$$ 
the universal subschemes, and set standardly $${\mathcal L}^{[d_v]} = \det {\mathbf R} p_{\star} \left ({\mathcal O}_{\mathcal Z} \otimes q^{\star} \mathcal L \right ), \, \,  \, {\mathcal L}^{[d_w]} = \det {\mathbf R} p_{\star} \left ({\mathcal O}_{\mathcal W} \otimes q^{\star} \mathcal L\right ).$$   From the product $${\mathcal X}^{[d_v]} \times_{{\mathcal K}_{\ell}} {\mathcal X}^{[d_w]} \times_{{\mathcal K}_{\ell}} {\mathcal X},$$ we calculate
\begin{eqnarray*}
\theta &= &\det  \left ( {\mathbf R p_{12}}_{\star} \left ( p_{13}^{\star} {\mathcal I}_{\mathcal Z} \otimes^{\mathbf L} p_{23}^{\star} {\mathcal I}_{\mathcal W} \otimes p_3^{\star} {\mathcal L} \right )\right ) ^{-1} \\
&=& \det  \left ( {\mathbf R p_{12}}_{\star} \left ( p_{13}^{\star}({\mathcal O} -  {\mathcal O}_{\mathcal Z}) \otimes^{\mathbf L} p_{23}^{\star}({\mathcal O} -  {\mathcal O}_{\mathcal W}) \otimes p_3^{\star} {\mathcal L} \right )\right ) ^{-1} \\
& = & {\mathcal L}^{[d_v]} \boxtimes {\mathcal L}^{[d_w]} \otimes \pi^{\star}( \det \pi_{\star} {\mathcal L} )^{-1}  \otimes \det  {\mathbf R p_{12}}_{\star} \left ( p_{13}^{\star} {\mathcal O}_{\mathcal Z} \otimes^{\mathbf L} p_{23}^{\star} {\mathcal O}_{\mathcal W} \otimes p_3^{\star} {\mathcal L} \right ) \\
& = &  {\mathcal L}^{[d_v]} \boxtimes {\mathcal L}^{[d_w]} \otimes \pi^{\star}( \det \pi_{\star} {\mathcal L} )^{-1}.  
\end{eqnarray*}
On the third line, the last bundle is the determinant of a complex of sheaves supported on the codimension 2 locus of intersecting subschemes in $ {\mathcal X}^{[d_v]} \times_{{\mathcal K}_{\ell}} {\mathcal X}^{[d_w]}$ -- thus it is trivial. Lemma $5.1$ of \cite {EGL} implies that $$\pi_{\star} \mathcal L^{[d_v]}=\Lambda^{[d_v]} \pi_{\star} \mathcal L,\,\,\, \pi_{\star} \mathcal L^{[d_w]}=\Lambda^{[d_w]} \pi_{\star} \mathcal L.$$ The higher direct images of the line bundles $\mathcal L^{[d_v]}, \mathcal L^{[d_w]}$ vanish by Theorem $5.2.1$ of \cite {Sc}. We therefore finally have
$$\pi_{\star} \theta \simeq \Lambda^{d_v} (\pi_{\star} {\mathcal L}) \otimes \Lambda^{d_w} (\pi_{\star} {\mathcal L}) \otimes (\det \pi_{\star} {\mathcal L} )^{-1}\cong \mathbf W'\otimes \mathbf V'.$$
We set $$\mathbf W'=\pi_{\star} \mathcal L^{[d_v]},\,\, \mathbf V'=\pi_{\star} \mathcal L^{[d_w]}\otimes (\det \pi_{\star} \mathcal L)^{\vee}.$$ As before these sheaves are only defined up to reverse twistings by a line bundle from $\mathcal P_1$. The divisor $\theta$ induces the duality isomorphism $$\mathsf D': \mathbf W'^{\vee}\to \mathbf V'$$ over $\mathcal P_1$, which is a global version of \eqref{sds}. 

Section $3$ shows that the universal relative Fourier-Mukai transform induces a birational map $$\mathcal X^{[d_v]}\times_{\mathcal P_1^{\circ}}\mathcal X^{[d_w]}\dasharrow M[v]\times_{\mathcal P_1^{\circ}} M[w]$$ regular in codimension $1$ over each fiber, such that the divisors $\theta$ and $\Theta$ are precisely matched. Because of regularity in codimension $1$, the pushforward sheaves  $\pi_{\star} \theta$ and $R^{0}\pi_{\star} \Theta$ coincide. Therefore $$\mathbf W'\otimes \mathbf V'\cong \mathcal H^0(\mathbf W)\otimes \mathcal H^0(\mathbf V)$$ over $\mathcal P_1^\circ$. We can furthermore align the line bundle twists inherent in the definition of $\mathbf W, \mathbf V, \mathbf W', \mathbf V'$ so that $$\mathcal H^0(\mathsf D)=\mathsf D'$$ over this locus. 
 We thus extended the Verlinde sheaves from $\mathcal P_1^{\circ}\hookrightarrow \mathcal P_1.$
 
 The resolution of the following query will however be of much greater interest.

\begin{question} Is it possible to extend $\mathbf W, \mathbf V$ from $$\mathcal K_{\ell}^{\circ}\hookrightarrow\mathcal K_{\ell}$$ in such a fashion that $$c_1(\mathbf W)=-c_1(\mathbf V)?$$
\end {question} 
\noindent Combined with the results of \cite {generic}, this would establish the strange duality conjecture over the entire locus where there is no higher cohomology, since the Baily-Borel compactification of $\mathcal K_{\ell}$ has one dimensional boundary. It would be interesting to investigate whether $\mathsf D$ is in fact a quasi-isomorphism between the complexes $\mathbf W^{\vee}$ and $\mathbf V$. 

Regarding the canonical line bundle $\Theta$, it is also natural to wonder 

\begin {question} Is the Chern character $\text {ch} ( \mathbf R\pi_{\star} \Theta)$ in the ring generated by the Hodge class $\lambda=-c_1(R^2\pi_{\star}\mathcal O_{\mathcal X^{\circ}})$  studied in \cite {kvg}?
\end {question}

\addtocounter{section}{+1}
\subsection {Acknowledgements} The authors were supported by NSF grants DMS 1001604 and DMS 1001486 and by the Sloan Foundation. Correspondence with G. van der Geer on related topics is gratefully acknowledged.

\end{document}